\pdfoutput=1
\documentclass[11pt]{article}
\usepackage{amsmath}
\usepackage{amssymb}
\usepackage{amsthm}
\usepackage{fancyhdr}
\usepackage[backend=biber, style=numeric]{biblatex}
\usepackage{multicol}
\newtheorem{theorem}{Theorem}[section]
\newtheorem{definition}[theorem]{Definition}
\newtheorem{lemma}[theorem]{Lemma}
\newtheorem{conjecture}[theorem]{Conjecture}
\newtheorem{corollary}[theorem]{Corollary}
\begin{document}
	\title{A conjecture on the prime factorization of $n!+1$}
	\author{William Gerst}
	\date{September 13, 2018}
	\maketitle

\begin{abstract}
	In this paper, we state a conjecture on the prime factorization of numbers of the form $n!+1$, explore its implications, and compare it with empirical evidence and established results based on the $abc$ conjecture.
\end{abstract}

\section{Introduction}
	Brocard's problem asks for integer solutions to the equation $n!+1=m^2$, and it is believed that there are only three such solutions. Known as Brown numbers, they are $(4,5)$, $(5,11)$, and $(7,71)$. In this paper, we will form a stronger conjecture on the factorization of numbers of the form $n!+1$, compare it to empirical results, and investigate the implications it would have regarding both the Brown numbers and another unsolved problem in number theory, the (in)finitude of the Wilson primes.

	To begin, we introduce two important functions pertaining to integer factorization, the sigma and omega functions.
	
\begin{definition}[Sigma function]
	Let $n\in\mathbb N$. The sigma function, or ``number-of-divisors'' function, denoted by $\sigma_0(n)$, returns the number of distinct divisors of $n$.	
\end{definition}

\begin{definition}[Omega function]
	Let $n\in\mathbb N$. The omega function, denoted by $\omega(n)$, returns the number of distinct prime divisors of $n$.
\end{definition}

	Note the difference between the functions $\sigma_0(n)$ and $\omega(n)$: the former counts the number of distinct divisors of $n$, while the latter counts only the number of distinct \textit{prime} divisors of $n$.

	Now that we have defined these functions, it is critical to note some of their most important characteristics.
	
\begin{lemma}
	Let $n\in\mathbb N$. Then $\sigma_0(n)$ is odd if and only if $n$ is a perfect square of an integer.
\end{lemma}
	
\begin{proof}
	Given a positive integer $n$ and a factor $d$, the number $n/d$ must also be a factor of $n$, by definition. Factors almost always come in pairs of the form $d$ and $n/d$, so the number of divisors will be even. However, if and only if there exists a $d$ such that $d=n/d$, then the two factors are no longer distinct, and so it is only counted once, thus changing the parity of the divisor count. Hence, the divisor count of any integer $n$ will be odd if and only if it is a perfect square.
\end{proof}

	For any $n\in\mathbb N$, we can calculate $\sigma_0(n)$ based on the prime factorization of $n$. Note that, by the fundamental theorem of arithmetic, the prime factorization of $n$ can be written
	
\begin{equation}
	n=\prod_{i=1}^{\omega(n)}(p_i)^{a_i}=(p_1)^{a_1}\cdot (p_2)^{a_2}\cdot (p_3)^{a_3}\cdots (p_{\omega(n)})^{a_{\omega(n)}},
\end{equation}
	where each $p_i$ represents a distinct prime and $a_i$ is the greatest power of the corresponding prime which divides $n$. Let $a_i$ be known as the \textbf{multiplicity} of the factor $p_i$, and note that $a_i>0$ in order for $p_i$ to divide $n$.
	
\begin{lemma}
	Let $n\in\mathbb N$. Then
	\begin{equation}
		\sigma_0(n)=\prod_{i=1}^{\omega(n)}(a_i+1).
	\end{equation}
\end{lemma}

\begin{proof}
	Using (1), $n=(p_1)^{a_1}\cdot (p_2)^{a_2}\cdot (p_3)^{a_3}\cdots (p_{\omega(n)})^{a_{\omega(n)}}$ for primes $p_i$ and multiplicities $a_i$. It is possible to construct a new number, \[k=\prod_{i=1}^{\omega(n)}(p_i)^{b_i}=(p_1)^{b_1}\cdot (p_2)^{b_2}\cdot (p_3)^{b_3}\cdots (p_{\omega(n)})^{b_{\omega(n)}},\] where $0\leq b_i\leq a_i$, which implies that $k$ is a divisor of $n$. For each prime factor $p_i$ where $1\leq i\leq\omega(n)$, there are $a_i+1$ choices of exponent, with each $b_i$ independent and each combination of choices producing a unique $k$, by the fundamental theorem of arithmetic. Therefore, the total number of distinct divisors of $n$ is given by \[\sigma_0(n)=\prod_{i=1}^{\omega(n)}(a_i+1).\]
\end{proof}

	There is one more definition that must be noted, and that is the concept of a square-free number.

\begin{definition}[Square-free numbers]
	A number $n\in\mathbb N$ is square-free if and only if there are no perfect squares which divide $n$ other than $1$.
\end{definition}

\section{The conjecture}
	Now that we have defined the concepts of square-free numbers and factor multiplicity, it is appropriate to introduce the conjecture. Henceforth, we will use $S$ to denote the set $\{4,5,7,12,23,229,562\}$.
	
\begin{conjecture}
	Let $n\in\mathbb N\setminus S$. Then $n!+1$ is square-free.
\end{conjecture}

	Some members of the excluded set $S$ should appear familiar. The first three, $4$, $5$, and $7$, are the only currently known solutions to Brocard's problem, mentioned earlier. In fact, each of these seven numbers has something in common which leads to their inclusion here: the prime factorization of $n!+1$ for any $n\in S$ produces a prime factor with multiplicity of at least $2$, and these are the only numbers of the required form currently known to have this property. As such, they present the only would-be counterexamples to the conjecture, as demonstrated below.

\begin{lemma}
	Let $n\in\mathbb N$. By (1), for primes $p$ and multiplicities $a$, \[n=\prod_{i=1}^{\omega(n)}(p_i)^{a_i}.\] If and only if $n$ is square-free, then $a_i=1$ for all $i$.
\end{lemma}

\begin{proof}
	If $n$ is square-free, then no perfect squares divide $n$. Assume for the sake of contradiction that for some $i$, $a_i>1$. Then, the perfect square $p_i^2$ divides the right-hand side of the above equation, as $(p_i)^2\mid (p_i)^{a_i}$ for $a_i>1$. Hence, we have found a perfect square that divides a square-free number, which is a contradiction. Therefore, there must be no $i$ such that $a_i>1$. By definition, $a_i>0$, so $n$ being square-free implies $a_i=1$ for all $i$.
	
	It remains to show that $a_i=1$ for all $i$ implies that $n$ is square-free. Assume for the sake of contradiction that $n$ is not square-free, and let $k\neq1$ be a perfect square that divides $n$. By definition, $k$ is twice-divisible by an integer square root, which itself is divisible by at least one prime. Then, $k$ is twice-divisible by that prime, and $n$ is as well. However, for each prime $p$ that divides $n$, $n$ is divisible by $p$ only once, a contradiction. Therefore, $a_i=1$ for all $i$ implies that $n$ is square-free.
\end{proof}

	Now that the multiplicity of each prime factor of a square-free number is known to be exactly $1$, it is possible to state one of the more important consequences.

\begin{theorem}
	Let $n\in\mathbb N$. If and only if $n$ is square-free, then
	\[\sigma_0(n)=2^{\omega(n)}.\]
\end{theorem}

\begin{proof}
	By equation (2) and lemma (2.2), if $n$ is square-free, then \[\sigma_0(n)=\prod_{i=1}^{\omega(n)}(a_i+1)=\prod_{i=1}^{\omega(n)}(1+1)=\prod_{i=1}^{\omega(n)}(2)=2^{\omega(n)}.\]
	
	It remains to be proven that $\sigma_0(n)=2^{\omega(n)}$ implies that $n$ is square-free. (The steps above cannot simply be reversed because $a_i=1$ implies the second equality in the chain, but this implication is not bidirectional. For instance, if $\omega(n)=2$, $a_1=0$, and $a_2=3$, the products would be equal but $a_i=1$ would not hold.) Note that $a_i>0$ by definition, so $(a_i+1)>1$. Again, using equation (2), \[2^{\omega(n)}=\sigma_0(n)=\prod_{i=1}^{\omega(n)}(a_i+1).\] As $2$ is prime, the only way to form $2^{\omega(n)}$ from a product of $\omega(n)$ integers greater than $1$ is as the product of $\omega(n)$ copies of $2$. The proof of this fact is trivial. Hence, it must be true that each multiplicand of $(a_i+1)$ in the product is equal to $2$, so $a_i=1$ for all $i$. Therefore, by lemma (2.2), $n$ must be square-free.
\end{proof}

\begin{corollary}
	Let $n\in\mathbb N$. If and only if $n!+1$ is square-free, then
	\[\sigma_0(n!+1)=2^{\omega(n!+1)}.\]
\end{corollary}

\begin{proof}
	This fact immediately follows from theorem (2.3), upon substituting $n!+1$ in place of $n$.
\end{proof}

	The significance of the above statement is astronomical, as it provides a way to quickly verify if $n!+1$ if square-free or not. If the conjecture (2.1) is true, then the equality $\sigma_0(n!+1)=2^{\omega(n!+1)}$ holds for all natural $n$ not contained in the excluded set $S$. We will now present numerical evidence which exemplifies this behavior and supports the claim set forth.

\section{Numerical evidence}
	There is very strong empirical evidence to support the results suggested in (2.4) for numbers of the form $n!+1$.
	
\begin{multicols}{2}
$\begin{array}{c|c|c}
	n & \sigma_0(n!+1) & 2^{\omega(n!+1)} \\
	\hline
	1 & 2 & 2 \\
	2 & 2 & 2 \\
	3 & 2 & 2 \\
	\mathbf{4} & \mathbf{3} & \mathbf{2} \\
	\mathbf{5} & \mathbf{3} & \mathbf{2} \\
	6 & 4 & 4 \\
	\mathbf{7} & \mathbf{3} & \mathbf{2} \\
	8 & 4 & 4 \\
	9 & 8 & 8 \\
	10 & 4 & 4 \\
	11 & 2 & 2 \\
	\mathbf{12} & \mathbf{6} & \mathbf{4} \\
	13 & 4 & 4 \\
	14 & 4 & 4 \\
	15 & 8 & 8 \\
	16 & 4 & 4 \\
	17 & 2 & 2 \\
	18 & 6 & 4 \\
	19 & 4 & 4 \\
	20 & 4 & 4 \\
\end{array}$
$\begin{array}{c|c|c}
	n & \sigma_0(n!+1) & 2^{\omega(n!+1)} \\
	\hline
	21 & 8 & 8 \\
	22 & 8 & 8 \\
	\mathbf{23} & \mathbf{12} & \mathbf{8} \\
	24 & 4 & 4 \\
	25 & 4 & 4 \\
	26 & 4 & 4 \\
	27 & 2 & 2 \\
	28 & 4 & 4 \\
	29 & 8 & 8 \\
	30 & 32 & 32 \\
	31 & 16 & 16 \\
	32 & 16 & 16 \\
	33 & 32 & 32 \\
	34 & 4 & 4 \\
	35 & 32 & 32 \\
	36 & 64 & 64 \\
	37 & 2 & 2 \\
	38 & 4 & 4 \\
	39 & 16 & 16 \\
	40 & 128 & 128 \\
\end{array}$
\end{multicols}
	
	In the table shown above, emphasis has been added to every row where $\sigma_0(n!+1)\neq2^{\omega(n!+1)}$, indicating the first five members of the excluded set $S$. It is confirmed in \cite{asahi} and \cite{OEIS} that this pattern holds for every $n$ between $41$ and $100$, as well. The next case known to break from the pattern is $n=229$, as $(229!+1)$ has a repeated prime factor of $613$ \cite{fact}. After that, the only other instance currently known is $562$, the final member of the excluded set $S$, where $(562!+1)$ is divisible by $563^2$ \cite{wilson}.
	
	Furthermore, it has been shown that the $abc$ conjecture implies that $n!+1=m$ has finitely many solutions for powerful $m$ \cite{power}, where powerful numbers are defined as those for which every prime divisor has a multiplicity greater than $1$. This implies there are infinitely many non-powerful $m$, where each is only once-divisible by \textit{at least one} of their prime divisors. The conjecture (2.1), on the other hand, postulates that there are infinite $m$ that are only once-divisible by \textit{all} of their prime divisors, and explicitly notes each of the finite exceptions to this rule.
	
\section{Further implications}

	It is important to return to Brocard's problem and the Brown numbers to understand the full implications of the conjecture at hand.
	
\begin{corollary}[Brocard's problem]
	If (2.1) is true, the only solutions to $n!+1=m^2$ for positive integers $(n,m)$ are $(4,5)$, $(5,11)$, and $(7,71)$.
\end{corollary}

\begin{proof}
	The conjecture (2.1) asserts that $n!+1$ is square-free unless $n$ is in the set $S$, and it is clear that a square-free number is, by definition, not a perfect square. Therefore, there are no perfect squares $n!+1$ when $n\notin S$. However, when $n\in S$, the only integers $n$ which satisfy the equation are $4$, $5$, and $7$. Hence, these three are the only solutions to the equation.
\end{proof}

	This conjecture implies that there are no more Brown numbers than the three pairs already discovered, as has been shown to hold for at least $n\leq10^9$ \cite{brocard}. In general, any result which shows that $\sigma_0(n!+1)$ is even for $n\in\mathbb N\setminus S$, by (1.3), would have this same implication. However, the conjecture we introduce in this paper is strictly stronger than that claim because it does not only assert that $\sigma_0(n!+1)$ is even, but also that it is a power of $2$.

	There is another significant consequence of the conjecture (2.1), and this one concerns the Wilson primes. A Wilson prime is a prime $p$ such that $p^2$ divides $(p-1)!+1$. There are believed to be infinitely many such primes, but $5$, $13$, and $563$ are known to be the only such numbers below $2\cdot10^{13}$ \cite{wilson}.
		
\begin{corollary}[Finitude of the Wilson primes]
	If (2.1) is true, the only Wilson primes are $5$, $13$, and $563$.
\end{corollary}

\begin{proof}
	According to conjecture (2.1), $n!+1$ is square-free if $n\notin S$. By definition, the square of a Wilson prime $p$ divides $(p-1)!+1$, so let $n=p-1$ and it is clear that $(p-1)\in S$ must be true in order for the conjecture to hold. It is easily verified that the only primes $p$ for which $(p-1)\in S$ are $5$, $13$, and $563$, and these are all Wilson primes. Hence, these three are the only Wilson primes that exist.
\end{proof}

	With the above two results, it is clear that a proof of this conjecture would directly and simultaneously imply the existence of only three pairs of Brown numbers and three Wilson primes, thus at once resolving two major unsolved problems in the field of number theory.

\printbibliography

\end{document}